\documentclass{amsart}
\newtheorem{statement}{\bf Statement}
\newtheorem{lemma}{\bf Lemma}
\newtheorem{corollary}{\bf Corollary}
\begin{document}
\title{A special Lipschitz map $f: \mathcal{C} \to H$ for the
Hadamard
space $H$. Sketch of the construction}
\author{P. N. Ivanshin}
\address{N.G. Chebotarev RIMM, KSU, Kazan, Russia}
\email{Pyotr.Ivanshin@ksu.ru}
\maketitle
The work on this subject was motivated by \cite{LPS}.
So let $H$ be an $n$-dimensional space of nonpositive curvature \cite{Bus}.  The problem is to find a Lipschitz function $f: \mathcal{C} \to H$ such that $f(\{x\})=\{x\}$ for any point $x \in H$. Assume that $H$ is $1$-connected and each two points of it can be connected with exactly one geodesic segment.
The construction is inductive one. The construction for the $1$-dimensional
space $H$ is evident.

1) Consider a point $p \in \partial_{\infty}H$ \cite{KL}.
Consider a limit $d_p(x, y) = \lim\limits_{t \to \infty} d(x(t),
y(t))$ for any pair of points $x, y$ on a horosphere $O_p^t$. Each
set $C \subset \mathcal{C}$ has a unique projection on $O_p$.
Consider a problem for the limit space $\lim\limits_{t \to
\infty}O_p(t), d_p$. Then take a first horosphere $O_1$
intersecting $C$ and take an image of the limit point on it.

The only natural obstacle is the existence of the subsets on
projections that shrink to a point in the limit. Note that any
such subset is a closed subset of $O_1$.

2a) If the set $C'$ shrinks as $t \to \infty$ then its preimage
with respect to the projection onto the $O_1$ is convex set and
has only one intersection point $p(C')$ with the appropriate
horosphere $O_{C'}$.

2b) If it does not then the map $d_p(x, y)$ is a correctly defined
metric. and there could be applied an induction step.

3) For the set of points $p_1, \ldots p_n , \ldots \subset O_1$ we
apply a center mass construction from \cite{I}. The only thing
left to show is that this center depends on points with mass in
Lipshitz way. To do this one must consider the convex set $C$ of
the fixed diameter $D$ and mass $M(C)$. The problem now can be
reduced to the case of the $0$ curvature and the case of strictly
negative one.

\begin{lemma}
Let $p \in \partial_{\infty} H$ If there exists a geodesic segment
on a horosphere $O_p$ then the distance between any two points on
this segment does not vanish for $t \to p$.
\end{lemma}

$\bullet$ In case an horosphere possess a geodesic segment $T(x,
y)$ there exists a Lambert rectangle of the type $\geq \pi/2$ then
by assumption it must be $\pi/2$ or equivalently the geodesic
lines do not converge one to another infinitely close, namely the
distance between them equals length $d(x, y)$. $\triangleright$

Center of mass construction.
\begin{lemma}
Consider a set $X=\bigcup\limits_{i=1}^{n}\{x_i\} \subset H$. The
diameter of the convex hull $\mathrm{diam}(Co(x_1, \ldots,
x_n))=\max\{d(x_i, x_j)| x_i, x_j \in X\}$.
\end{lemma}

\begin{corollary}
Suppose for two sets $X= \bigcup\limits_{i=1}^{n}\{x_i\}$ and $Y=\bigcup\limits_{i=1}^{n}\{y_i\}$, $x_i, y_i \in H$ $\exists a>0$ the following: $\forall i, j \in \{1, \ldots, n\}$ $\exists k, l \in \{1, \ldots, n\}$ $d(y_i, y_j) \leq a d (x_k, x_l)$ then $a\mathrm{diam}(Co(x_1, \ldots, x_n)) \geq \mathrm{diam}(Co(y_1, \ldots, y_n))$.
\end{corollary}

Now define a center of the mass of the set of points with mass as
follows: for two points $x_1, m_1$, $x_2, m_2$ this center is the
point on the geodesic segment connecting $x_1$ to $x_2$ dividing
the length of a segment in the proportion $m_1/m_2$. For $n$
points one must consider first the center of the mass of the $n-1$
points $(x_2, m_2), \ldots, (x_n, m_n)$ and then the first step of
the construction for these points $(x_1^{'},
\frac{\sum\limits_{i=2}^n m_i}{n-1}) \ldots, (x_n^{'},
\frac{\sum\limits_{i=1}^{n-1} m_i}{n-1})$. Consider then the same
construction for a set of centers $X_1$ with the given masses.
Then repeat the procedure.

\begin{statement}
The center of the set of points with masses $X \subset H$ is
correctly defined, i.e. it exists and is unique.
\end{statement}

Thus we can construct relatively good function for centers of mass
which occur in nonsingular points. In the singular points this
function will have discontinuities. To avoid this we must deform
our function as follows: for all $x \in O_{\infty}$ is a center
point for a set of the zero curvature such that the center of the
mass of the whole set of points lies in the singular point $y \in
O_{\infty}$ we must consider a function locally constant in the
neighbourhood of this point, so that the center of mass stays in
the point $y$. We need this to define a point on the fixed
horosphere $O_1$ so that if the center depending on $x$ equals $y$
then the target point on $O_1$ lies in the subsequent shrinking
set in the neighbourhood of the center point for this set.

The problem now reduces to the proof of the Lipshitz structure of
a center mass function. To do this we must consider two different
cases: 1) shift of one of the points of the set $X$; 2) change of
the mass of one of the points of $X$. The proof is again by
induction. It is trivial in case $|X|=1$. Then for $n$ points one
has on the first step of the construction a set of points. Note
that these points depend on the shift of one of them in Lipschitz
way. Since this holds true on each subsequent construction step we
can consider a converging series each element of which can be
majored by an element from the standard Euclidean space
$\mathbb{R}^n$ since the distance function for two points on distinct geodesic lines is convex \cite{Bus}.

\end{document}